\documentclass[12pt,a4paper]{article}
\usepackage{amsmath,amssymb,amsthm}
\usepackage{fullpage}

\newtheorem{thm}{Theorem}

\newcommand{\dis}{\displaystyle}

\newtheorem{lemma}{Lemma}
\newtheorem{Def}{Definition}
\numberwithin{equation}{subsection}
\begin{document}

\title{ON THE ARCHIMEDEAN OR SEMIREGULAR POLYHEDRA}
\author{Mark B. Villarino\\
Depto.\ de Matem\'atica, Universidad de Costa Rica,\\
2060 San Jos\'e, Costa Rica}
\date{May 11, 2005}

\maketitle

 \begin{abstract}
 We prove that there are thirteen \textsc{Archimedean}/semiregular polyhedra by using \textsc{Euler}'s polyhedral formula.
 \end{abstract}
\tableofcontents
\section{Introduction}
\subsection{Regular Polyhedra}

A \textbf{\emph{polyhedron}} may be intuitively conceived as a ``solid figure" bounded by plane faces and straight line edges so arranged that every edge joins exactly \emph{two} (no more, no less) vertices and is a common side of two faces.

A polyhedron is \textbf{\emph{regular}} if all its \textbf{\emph{faces}} are regular polygons (with the same number of sides) and all its \textbf{\emph{vertices}} are regular polyhedral angles; that is to say, all the face angles at every vertex are congruent and all the dihedral angles are congruent.  An immediate consequence of the definition is that all the faces of the polyhedron are congruent.

\textbf{\emph{There are FIVE such regular convex polyhedra}}, a fact known since \textsc{Plato}'s time, at least, and all of Book XIII of \textsc{Euclid} is devoted to proving it, as well as showing how to construct them:  the \textbf{\emph{tetrahedron}}, the \textbf{\emph{cube}}, the \textbf{\emph{octahedron}}, the \textbf{\emph{dodecahedron}}, and the \textbf{\emph{icosahedron.}}

The following table summarizes the basic data on the regular polyhedra.

\large
\begin{center}
\fbox{\textbf{Regular Polyhedra}}
\end{center}
\normalsize
\begin{center}
\begin{tabular}{|c|c|c|c|c|c|}
  \hline
   &\textbf{POLYGONS}  &  &  &  & \textbf{NUMBER OF}  \\
 \textbf {NAME } &\textbf{ FORMING} &\textbf{VERTICES} & \textbf{EDGES}& \textbf{FACES}&\textbf{FACES AT EACH}\\
 &\textbf{FACES} & & & & \textbf{VERTEX}\\
  \hline
 \textbf{Tetrahedron} & \textbf{Triangles} & $4$ & $6$ & $4$ & $3$\\ \hline
 \textbf{Octahedron} & \textbf{Triangles} & $6$ & $12$ & $8$ & $4$\\ \hline
 \textbf{Icosahedron} & \textbf{Triangles} & $12$ & $30$ & $20$ & $5$\\ \hline
 \textbf{Cube} & \textbf{Squares} & $8$ & $12$ & $6$ & $3$\\ \hline
 \textbf{Dodecahedron} & \textbf{Pentagons} & $20$ & $30$ & $12$ & $3$\\ \hline
\end{tabular}
\end{center}
\subsection{Archimedean/semiregular polyhedra}

It's reasonable to ask what happens if we \emph{forgo} some of the conditions for regularity.  \textsc{Archimedes} \cite{Arch} investigated the polyhedra that arise if we \emph{retain} the condition that the faces have to be regular polygons, but \emph{replace} the regularity of the polyhedral angles at each vertex by the \emph{weaker} condition that they all be congruent (see \textsc{Lines} \cite{Lines}).  Such solids are called \textbf{\emph{Archimedean}} or \textbf{\emph{semiregular}} polyhedra.

\begin{thm}(\textbf{Archimedes' Theorem})
There are \textbf{THIRTEEN} semiregular polyhedra as well as two infinite groups: the prisms and the antiprisms.
\end{thm}
  The following table summarizes the data on the thirteen semiregular polyhedra. The labels are self-explanatory except  for the \emph{C \& R symbol} \cite{CR}: $\mathbf{a^{b}.c^{d}}...$ \emph{means $\mathbf{b}$ regular $\mathbf{a}$-gons, $\mathbf{d}$ regular $\mathbf{c}$-gons, . . . meet at a vertex.}  Moreover the number of regular $k$-gonal facets is denoted by $F_{k}$.  

Thus, for example, the \emph{cuboctahedron} has $12$ vertices, $24$ edges, $14$ facets, of which $8$ are equilateral triangles and $6$ are squares.  Moreover, at each vertex one has a triangle, a square, a triangle, a square, in that cyclic order.

\large
\begin{center}
\fbox{\textbf{Archimedean Polyhedra}}
\end{center}
\normalsize
\begin{center}
\begin{tabular}{|c|c|c|c|c|c|c|c|c|c|c|}
  \hline
   & & & & & & & & & & \textbf{C \& R}\\
  \textbf{NAME} &\textbf{V}  & \textbf{E} & \textbf{F} &$\mathbf{F_{3}}$  & $\mathbf{F_{4}}$ &$\mathbf{F_{5}}$&$\mathbf{F_{6}}$ &$\mathbf{F_{8}}$ & $\mathbf{F_{10}}$&   \\
  & & & & & & & & & & \textbf{Symbol}\\
 \hline
\textbf{cuboctahedron} &12 &24 &14 &8 & 6& & & & & $(3.4)^{2}$\\
 \hline
 \textbf{great rhombicosidodecahedron} &120 &180 &62 & & 30& &20 & &12 & $4.6.10$\\
 \hline
 \textbf{great rhombicuboctahedron} &48 &72 &26 & & 12& &8 & 6& & $4.6.8$\\
 \hline
 \textbf{icosidodecahedron} &30 &60 &32 &20 & &12 & & & & $(3.5)^{2}$\\
 \hline
 \textbf{small rhombicosidodecahedron} &60 &120 &62 &20 & 30&12 & & & & $(3.4.5.4)$\\
 \hline
 \textbf{small rhombicuboctahedron} &24 &48 &26 &8 & 18& & & & & $3.4^{3}$\\
 \hline
 \textbf{snub cube}&24 &60 &38 &32 & 6& & & & & $3^{4}.4$\\
 \hline
 \textbf{snub dodecahedron} &60 &150 &92 &80 & &12 & & & & $3^{4}.5$\\
 \hline
 \textbf{truncated cube} &24 &36 &14 &8 & & & &6 & & $3.8^{2}$\\
 \hline
 \textbf{truncated dodecahedron} &60 &90 &32 &20 & & & & &12 & $3.10^{2}$\\
 \hline
 \textbf{truncated icosahedron} &60 &90 &32 & & &12 &20 & & & $4.6.10$\\
 \hline
 \textbf{truncated octahedron} &24 &36 &14 & & 6& &8 & & & $4.6^{2}$\\
 \hline
 \textbf{truncated tetrahedron} &12 &18 &8 &4 & & &4 & & & $3.6^{2}$\\
 \hline
\end{tabular}
\end{center}  
\section{Proof techniques}
\subsection{Euclid's proof for regular polyhedra}

\textsc{Euclid}'s proof (Proposition XVIII, Book XIII) is based on the \textbf{\emph{polyhedral angle inequality}}:  \emph{the sum of the face angles at a vertex cannot exceed }$2\pi$, as well as on the fact that \emph{the internal angle of a regular $p$-gon is} $\pi-\frac{2\pi}{p}.$

Thus, if $q$ faces meet at each vertex
\begin{align}
&\Rightarrow q\left(\pi-\frac{2\pi}{p}\right)  <  2\pi \\
 &\Rightarrow (p-2)(q-2)  <  4 \\
& \Rightarrow (p,q)=(3,3),\ (4,3),\ (3,4),\ (5,3),\ (3,5)
\end{align}
which give the tetrahedron, cube, octahedron, dodecahedron, and icosahedron respectively.

Of course the key step is to obtain (2.1.2).  \textsc{Euclid} does it by (2.1.1) which expresses a \emph{metrical} relation among angle measures.

One presumes that \textsc{Archimedes} applied more complex versions of (2.1.1) and (2.1.2)
to prove that the semiregular solids are those thirteen already listed.  Unfortunately, his treatise was lost over two thousand years ago!

\subsection{Euler's polyhedral formula for regular polyhedra}

Almost the same amount of time passed before somebody came up with an entirely new proof of (2.1.2), and therefore of (2.1.3).  In 1752 \textsc{Euler}, \cite{Euler}, published his famous \textbf{\emph{polyhedral formula:}}

\begin{equation}
 \fbox{$\dis V-E+F=2$}
\end{equation}
in which $V:=$ the number of vertices of the polyhedron, $E:=$ the number of edges, and $F:=$ the number of faces.  This formula is valid for any polyhedron that is homeomorphic to a sphere.

The proof of (2.1.2) using (2.2.1) goes as follows.  If $q$ $p$-gons meet at each vertex,
\begin{align}
 & \Rightarrow pF=2E=qV  \\
&\Rightarrow E=\frac{qV}{2}, \ F=\frac{qV}{p} 
\end{align}Substituting (2.2.3) into (2.2.1), \begin{align*}
&\Rightarrow V-\frac{qV}{2}+\frac{qV}{p}=2 \\
&\Rightarrow 2pV-qpV+2qV=4p\\
&\Rightarrow V=\frac{4p}{2p-qp+2q}\\
&\Rightarrow 2p-qp+2q>0\\
&\Rightarrow (p-2)(q-2)<4 
\end{align*}which is (2.1.2).

This second proof proves much more.  We have found \emph{all regular maps} (graphs, networks) on the surface of a sphere whatever the boundaries may be, without \emph{any} assumptions in regard to they're being circles or skew curves.  Moreover the exact shape of the sphere is immaterial for our statements, which hold on a cube or any hemeomorph of the sphere.

This \emph{topological} proof of (2.1.2) is famous and can be found in numerous accessible sources, for example \textsc{Rademacher \& Toeplitz} \cite{R-T}.

\subsection{Proofs of Archimedes' theorem}

Euclidean-type \emph{metrical} proofs of \textsc{Archimedes}' theorem are available in the literature (see \textsc{Cromwell} \cite{Crom} and \textsc{Lines} \cite{Lines}) and take their origin in a proof due to \textsc{Kepler} \cite{Kepler}.

They use the polyhedral angle inequality to prove:\begin{itemize}
  \item at most \textbf{\emph{three}} different kinds of face polygons can appear around any solid angle;
  \item three polygons of different kinds \emph{cannot} form a solid angle if any of them has an \textbf{\emph{odd}} number of sides
  \end{itemize}One then exhaustively examines all possible cases.

The situation is quite different with respect to a \emph{topological} proof of \textsc{Archimedes}' theorem.  Indeed, after we had developed our proof, as presented in this paper, we were able to find only one reference: \textsc{T.R.S. Walsh} \cite{Walsh} in 1972.

His proof, too, is based exclusively on \textsc{Euler}'s polyhedral formula, and so there are overlaps with ours.  However, our proof is quite different, both in arrangement and details, and in purpose.  The pedagogical side is insisted upon in our proof so as to make it as elementary and self-contained as possible for as wide an audience as possible.
We comment further on the structure of this proof after the proof of \textbf{Lemma 3..}

\section{Three lemmas}

For any polyhedron we define:\begin{itemize}
  \item $V$:= total number of vertices;
  \item $V_{p}$:= total number of vertices incident with $p$ edges;
  \item $E$:= total number of edges;
  \item $F$:= total number of faces;
  \item $F_{p}$:= total number of $p$-gonal faces;\end{itemize}
Here, and from now on, \textbf{\emph{polyhedron}} means any map on the sphere for which \textsc{Euler}'s theorem holds.


\subsection{Lemma 1}

The following lemma is due to \textsc{Euler} \cite{Euler} and is well known.  We sketch the proof for completeness.

\begin{lemma}
The following relations are valid in any polyhedron:\begin{enumerate}
  \item $3F_{3}+2F_{4}+F_{5}=12+2V_{4}+4V_{5}+\cdots+F_{7}+2F_{8}+\cdots.$
  \item At least one face has to be a \textbf{triangle}, or a \textbf{quadrilateral}, or a \textbf{pentagon}; i.e., there is \textbf{no} polyhedron whose faces are all hexagons, or polygons with \textbf{six} or more sides.
\end{enumerate}
\end{lemma}

\textbf{Proof:}

For \textbf{1.} we note\begin{itemize}
  \item (i) $F_{3}+F_{4}+\cdots+F_{V-1}=F;$
  \item (ii) $3F_{3}+4F_{4}+\cdots+(V-1)F_{V-1}=2E;$
  \item (iii) $V_{3}+V_{4}+\cdots+F_{V-1}=V;$
  \item (iv) $3V_{3}+4V_{4}+\cdots+(F-1)V_{F-1}=2E.$
\end{itemize}Now multiply (i) by $6$, subtract (ii), and use (iii), (iv), and \textsc{Euler}'s formula.
\\

For \textbf{2.} observe that $F_{3}$, $F_{4}$, and $F_{5}$ cannot all be zero in \textbf{1.} at the same time.

\hfill$\Box$

\subsection{Lemma 2}

\begin{Def}
A polyhedron is called \textbf{\emph{Archimedean}} or \textbf{\emph{semiregular}} if the cyclic order of the degrees of the faces surrounding each vertex is the \textbf{\emph{same}} to within rotation and reflection. \cite{Walsh}
\end{Def}

\begin{lemma}
In any Archimedean polyhedron:\begin{enumerate}
  \item \begin{equation*}
\fbox{$\dis rV=2E$}
\end{equation*} where $r$ edges are incident at each vertex.
  \item \begin{equation*}
\fbox{$\dis \frac{pF_{q}}{q}=V$}
\end{equation*}where $q$ $p$-gons are incident at each vertex.
  \item \begin{equation*}
\fbox{$\dis V=\frac{2}{ 1-\frac{r}{2}+\frac{1}{p_{1}}+\frac{1}{p_{2}}+\cdots+\frac{1}{p_{r}}}$}
\end{equation*}where one $p_{1}$-gon,  one $p_{2}$-gon, $\cdots$, one $p_{r}$-gon all meet at one vertex and where the $p_{k}$ don't all have to be different.  
\end{enumerate}
\end{lemma}

\textbf{Proof:}

For \textbf{1.}, since there are $2$ vertices on any edge, the product $rV$ counts each edge twice, so $=2E.$
\\

For \textbf{2.}, $pF_{p}$ counts the total number of vertices \emph{once} if \emph{one} $p$-gon is incident at each vertex, \emph{twice} if \emph{twp} $p$-gons are incident there, $\cdots$, $q$ times if $q$ $p$-gons are incident at the vertex.  That is, $pF_{p}=qV.$
\\

For \textbf{3.}, solve \textbf{1.} for $E$, use (i) of the proof of \textbf{Lemma 1.1}, solve \textbf{2.} for $F_{p},$ substitute in \textsc{Euler}'s formula, solve for $V$, and write any fraction $$\frac{q}{p}=\underbrace{\frac{1}{p}+\frac{1}{p}+\cdots+\frac{1}{p}}_{q\  \rm times}.$$

\hfill$\Box$


\subsection{Lemma 3}

This lemma limits the number of candidate polygons surrounding each vertex.

\begin{lemma}
If $r$ edges are incident with each vertex of an Archimedean polyhedron then\begin{equation*}
\fbox{$\dis r\leqslant 5$}
\end{equation*}
\end{lemma}

\textbf{Proof:}

By \textbf{3.} of \textbf{Lemma 2.}\begin{align*}
&1-\frac{r}{2}+\frac{1}{p_{1}}+\cdots+\frac{1}{p_{r}}>0\\
&\Rightarrow \frac{1}{p_{1}}+\cdots+\frac{1}{p_{r}}>\frac{r-2}{2}
\end{align*}But, \begin{align*}\\
&p_{1}\geqslant 3, \ p_{2}\geqslant 3, \cdots, \ p_{r}\geqslant 3\\
&\Rightarrow \frac{1}{3}+\frac{1}{3}+\cdots+\frac{1}{3}\geqslant \frac{1}{p_{1}}+\cdots+\frac{1}{p_{r}}> \frac{r-2}{2}\\
&\Rightarrow\frac{r}{3}>\frac{r-2}{2}\\
&\Rightarrow r<6\\
&\Rightarrow r\leqslant 5.
\end{align*}

\hfill$\Box$

It is of interest to compare the method of proof, using \textsc{Euler}'s theorem, for the \textbf{\emph{regular}} polyhedra and the \textbf{\emph{Archimedean}} polyhedra.

In \textbf{\emph{both}} cases the essential step is to use the fact that \textbf{\emph{the denominator of the formula for the number of vertices, V, is positive:}}\begin{align*}
V&=\frac{2}{1-\frac{r}{2}+\frac{r}{p}}\ \ \ \ \ \ \ \ \ \ \ \ \ \ \ \ \ \ \ \ \ \ \ \ \ \ \ \ \ \ \rm regular\\
V&=\frac{2}{1-\frac{r}{2}+\frac{1}{p_{1}}+\frac{1}{p_{2}}+\cdots+\frac{1}{p_{r}}}\ \ \ \ \ \ \ \ \rm Archimedean
\end{align*}In the case of the \emph{regular} polyhedron the inequality $$1-\frac{r}{2}+\frac{r}{p}>0$$can be rearranged into the elegant inequality $$(p-2)(r-2)<4,$$which, as we saw before, leads to five solutions $(p,r).$

Unfortunately, in the case of the \emph{Archimedean} polyhedra the inequality $$1-\frac{r}{2}+\frac{1}{p_{1}}+\frac{1}{p_{2}}+\cdots+\frac{1}{p_{r}}>0$$apparently does \textbf{\emph{not}} lend itself to an algebraic rearrangement into a product, and so must be studied by \emph{an exhaustive enumeration of cases}. 

Nevertheless, it's worth emphasizing that the basic structure of the two arguments is the same at the core, although the elaboration of the cases in the Archimedean case demands some topological counting arguments that are not entirely trivial.


\section{Topological Proof of Archimedes' theorem}

By \textbf{Lemma 3} we have to consider three cases:

\begin{itemize}
  \item Case 1: Five faces meet at a vertex $r=5.$
  \item Case 2: Four faces meet at a vertex $r=4.$
  \item Case 3: Three faces meet at a vertex $r=3.$
\end{itemize}

\subsection{Case 1: five faces meet at a vertex: r=5}

By \textbf{Lemma 3.2},\begin{align*}&1-\frac{5}{2}+\frac{1}{p_{1}}+\frac{1}{p_{2}}+\frac{1}{p_{3}}+\frac{1}{p_{4}}+\frac{1}{p_{5}}=\frac{2}{V}>0\end{align*}\begin{align}
&\Rightarrow \frac{1}{p_{1}}+\frac{1}{p_{2}}+\frac{1}{p_{3}}+\frac{1}{p_{4}}+\frac{1}{p_{5}}-\frac{3}{2}>0
\end{align}By \textbf{Lemma 1.2}, at least one of $p_{1}$, $\cdots$, $p_{5}$ has to be $3$, $4$, or $5$.
\subsubsection{At least one face is a triangle: $p_{1}=3$}

Assuming $p_{1}=3$, \begin{align*}
&\Rightarrow \frac{1}{p_{2}}+\frac{1}{p_{3}}+\frac{1}{p_{4}}+\frac{1}{p_{5}}-\frac{3}{2}+\frac{1}{3}>0\\
&\Rightarrow \frac{1}{p_{2}}+\frac{1}{p_{3}}+\frac{1}{p_{4}}+\frac{1}{p_{5}}-\frac{7}{6}>0
\end{align*}Without loss of generality, we assume that:\begin{align*}
&3\leqslant p_{2}\leqslant p_{3}\leqslant p_{4}\leqslant p_{5}\\
&\Rightarrow \frac{1}{3}\geqslant\frac{1}{p_{2}}\geqslant \frac{1}{p_{3}}\geqslant\frac{1}{p_{4}}\geqslant\frac{1}{p_{5}}\\
&\Rightarrow \frac{1}{3}+\frac{1}{3}+\frac{1}{3}+\frac{1}{p_{5}}-\frac{7}{6}>0\\
&\Rightarrow \frac{1}{p_{5}}-\frac{1}{6}>0\\
&\Rightarrow p_{5}<6\\
&\Rightarrow p_{5}= 5, \ 4, \ 3\\
&\Rightarrow (p_{1},p_{2},p_{3},p_{4},p_{5})=(3,p_{2},p_{3},p_{4},5), (3,p_{2},p_{3},p_{4},4), (3,p_{2},p_{3},p_{4},3)
\end{align*}However, if we take \begin{align*}
&p_{2}\geqslant 3, p_{3}\geqslant 3, p_{4}\geqslant 4, p_{5}\geqslant 4\\
&\Rightarrow \frac{1}{p_{1}}+\frac{1}{p_{2}}+\frac{1}{p_{3}}+\frac{1}{p_{4}}+\frac{1}{p_{5}}\leqslant \frac{1}{3}+\frac{1}{3}+\frac{1}{3}+\frac{1}{4}+\frac{1}{4}=\frac{3}{2}
\end{align*}and this contradicts (4.1.1). Therefore we are left with only three quintuplets:
\begin{equation}
\fbox{$\dis (p_{1},p_{2},p_{3},p_{4},p_{5})=(3,3,3,3,5), (3,3,3,3,4), (3,3,3,3,3). $}
\end{equation} These correspond, respectively, to the \textbf{\emph{snub dodecahedron}}, the \textbf{\emph{snub cube}}, and the \textbf{\emph{icosahedron}}, a regular polyhedron.

\subsubsection{All faces have at least four sides: $p_{1}\geqslant 4$}

Again, we must assume \begin{align*}
&4\leqslant p_{1}\leqslant p_{2}\leqslant p_{3}\leqslant p_{4}\leqslant p_{5}\\
&\Rightarrow \frac{1}{4}+\frac{1}{4}+\frac{1}{4}+\frac{1}{4}+\frac{1}{p_{5}}-\frac{3}{2}>0\\
&\Rightarrow \frac{1}{p_{5}}-\frac{1}{2}>0\\
&\Rightarrow p_{5}<2 \ (\Rightarrow \Leftarrow) 
\end{align*}Since one of the cases $p_{1}=3$ or $p_{1}\geqslant 4$ must hold, and since they exhaust all possiblities with $r=5$, we are left with:
\begin{equation}
\fbox{$\begin{array}   {rcl}   (p_{1},p_{2},p_{3},p_{4},p_{5})&=&3^{4}.5\cdots \rm \mathbf{snub\  dodecahedron}\\
&=&3^{4}.4\cdots \rm \mathbf{snub cube}\\
&=&3^{5}\cdots \rm \mathbf{regular\  icosahedron}
\end{array}$}
\end{equation}
\subsection{Case 2: four faces meet at a vertex:  r=4}

By \textbf{Lemma 2.3}, \begin{align*}
&1-\frac{4}{2}+\frac{1}{p_{1}}+\frac{1}{p_{2}}+\frac{1}{p_{3}}+\frac{1}{p_{4}}>0\\
&\Rightarrow \frac{1}{p_{1}}+\frac{1}{p_{2}}+\frac{1}{p_{3}}+\frac{1}{p_{4}}-1>0.
\end{align*}Again, at least one of the $p_{k}$ must be $3$, $4$, or $5$.
\subsubsection{At least one face is a triangle: $p_{1}=3$}

We will write $p,q,r$ instead of $p_{1},p_{2},p_{3}.$ Thus the inequality becomes \begin{equation}\frac{1}{p}+\frac{1}{q}+\frac{1}{r}-\frac{2}{3}>0.\end{equation}

We examine a typical polyhedron: \begin{itemize}
  \item it must have a triangle at each vertex;
  \item there must be $4$ edges incident at each vetex;
  \item the vertices must all have the same configuration in the same order to within rotation and reflection.
\end{itemize}As we label the faces round each vertex of a triangle, say counterclockwise, from the lower right, we are compelled to conclude that \emph{no matter how we label the vertices} \textbf{\emph{at least two of the $\mathbf{p, \ q, \ r}$ must be equal.}} Putting $r=p$ in the inequality (4.2.1), we obtain \begin{align*}
 &\frac{2}{p}+\frac{1}{q}-\frac{2}{3}>0   \\
    &\Rightarrow (p-3)(2q-3)<9\\
    &\Rightarrow 1<2q-3<9, \ (2q-3)\ \ \rm odd\\
    &\Rightarrow 2q-3=3, \ 5,\  \ 7\\ 
    &\Rightarrow\begin{cases}
      & 2q-3=5 \  or \ 7, \Rightarrow p-3=0, \  \ 1 \Rightarrow p=3, \  \ 4 \ ;\text{otherwise}\\
      & 2q-3=3 \ \Rightarrow p-3=0, \ 1, \ 2 \ \Rightarrow p=3, \ 4, \ 5.
\end{cases} 
\end{align*}Therefore, we obtain\begin{center}\begin{tabular}{|c|c|c|c|c|c|c|c|}
\hline
  p &3  &3 &3 &4 &4 &4 &5 \\ 
  \hline
  q &3  &4 &5 &3 & 4& 5 &3\\ 
  \hline
\end{tabular}\end{center}Finally we observe that $2q-3\geqslant 9$ is permitted if $p-3=0.$

Therefore, we are left with:\begin{align*}
   \fbox{$ \begin{array}{lll}& (p,q)=(4,5)\Rightarrow (p_{1},p_{2},p_{3},p_{4}) = (3.4.5.4)\cdots\rm \ \mathbf{small} \ \mathbf{rombicosidodecahedron} \\
    &  (p,q)=(5,3)\Rightarrow (p_{1},p_{2},p_{3},p_{4}) = (3.5)^{2}\cdots\rm \ \mathbf{icosidodecahedron} \\
     &  (p,q)=(4,4)\Rightarrow (p_{1},p_{2},p_{3},p_{4}) = 3.4^{3}\cdots\rm \ \mathbf{small \ rhombicuboctahedron} \\
      &  (p,q)=(4,3)\Rightarrow (p_{1},p_{2},p_{3},p_{4}) = (3.4)^{2}\cdots\rm \mathbf{ cuboctahedron }\\
       &  (p,q)=(3,3)\Rightarrow (p_{1},p_{2},p_{3},p_{4}) = 3^{4}\cdots\rm \ \mathbf{regular \ octahedron} \\
        &  (p,q)=(3,m)\Rightarrow (p_{1},p_{2},p_{3},p_{4}) = 3^{3}.m \ (m\geqslant 4)\cdots\rm \ \mathbf{antiprism} \end{array}$}
\end{align*}
\subsubsection{All faces have at least four sides: $p_{1}\geqslant 4$}

We assume that\begin{align*}
&4\leqslant p_{1}\leqslant p_{2}\leqslant p_{3}\leqslant p_{4}   \\
    & \Rightarrow \frac{1}{4}+\frac{1}{4}+\frac{1}{4}+\frac{1}{p_{4}}>0\\
    &\Rightarrow p_{4}<4 \ ( But,\ p_{4}\geqslant 4 (\Rightarrow \Leftarrow)).
\end{align*}Therefore $p_{1}\geqslant 4$ can't happen.

There are no other cases with $r=4$.  

\subsection{Case 3: three faces meet at a vertes:  r=3}

By \textbf{Lemma 2.3},\begin{align*}
    &1-\frac{3}{2}+\frac{1}{p_{1}}+\frac{1}{p_{2}}+\frac{1}{p_{3}}>0\\
    & \Rightarrow  \frac{1}{p_{1}}+\frac{1}{p_{2}}+\frac{1}{p_{3}}-\frac{1}{2}>0
\end{align*}Since at least one of the $p_{k}$ must be equal to $3$, $4$, or $5$, we consider each case separately.

\subsubsection{At least one face is a triangle: $p_{1}=3$}

Then,\begin{align*}
    &  \frac{1}{p_{2}}+\frac{1}{p_{3}}-\frac{1}{6}>0 \\
\end{align*}

Looking at the configuration we see:\begin{itemize}
  \item each vertex has three edges incident to it
  \item two are the edges of a triangle and the third of a $p_{3}$-gonal face 
\end{itemize}  Labeling it we see that $$p_{2}=p_{3},$$
and therefore the above equality becomes \begin{align*}
    &\frac{2}{p_{3}}-\frac{1}{6}>0  \\
    & \Rightarrow p_{3}<12\\
    &3\leqslant p_{3}\leqslant 11. 
\end{align*}
\begin{lemma}
$p_{3}$ is \textbf{even} or $p_{3}=3.$
\end{lemma}
\begin{proof}
We look at the configuration with $p_{3}\geqslant 4.$ Since the vertices must all look alike, as we traverse counterclockwise (say) the $p_{3}$ vertices of a $p_{3}$-gonal face, we observe that the edges of the face fall into two groups:\begin{itemize}
  \item those that are the common edge of two $p_{3}$-gonal faces;
  \item those that are the common edge of a triangle and a $p_{3}$-gonal face. 
\end{itemize}Moreover, they occur in adjacent pairs, and finally, as we complete one circuit and return to our starting point, having started with a triangular edge, we end up with an edge common to two $p_{3}$-gonal faces.  Thus we traverse an \emph{integral number of \textbf{pairs} of sides} as we percorse the $p_{3}$-gonal face once, i.e., $p_{3}$ is \textbf{\emph{even}}.
\end{proof}

The only even numbers $p_{3}$ between $3$ and $11$ are $$p_{3}=4, \ 6, \ 8, \ 10.$$
Therefore we obtain

\begin{align*}
   \fbox{$ \begin{array}{lll}
    &  p_{3}=3\Rightarrow (p_{1},p_{2},p_{3}) = 3^{3}\cdots\rm \ \mathbf{regular \ tetrahedron} \\
     &  p_{3}=4\Rightarrow (p_{1},p_{2},p_{3}) = 3.4^{2}\cdots\rm \ \mathbf{triangular \ prism} \\
      &  p_{3}=6\Rightarrow (p_{1},p_{2},p_{3}) = 3.6^{2}\cdots\rm \ \mathbf{truncated \  tetrahedron }\\
       &  p_{3}=8\Rightarrow (p_{1},p_{2},p_{3}) =3.8^{2}\cdots\rm \ \mathbf{truncated \ cube} \\
        &  p_{3}=10\Rightarrow (p_{1},p_{2},p_{3}) = 3.10^{2} \ \cdots\rm \ \mathbf{truncated \ dodecahedron} \end{array}$}
\end{align*}


\subsubsection{All faces have at least four sides and one exactly four sides: $p_{1}=4\leqslant p_{2}\leqslant p_{3}.$}

Then, \begin{align*}
    &\frac{1}{4}+\frac{1}{p_{2}}+\frac{1}{p_{3}}-\frac{1}{2}>0   \\
    & \Rightarrow (p_{2}-4)(p_{3}-4)<16. 
\end{align*}The same sort of configuration argument shows \emph{that $p_{2}$ and $p_{3}$ are} 
\textbf{\emph{even}}.  Thus,\begin{align*}
    & p_{2}=2a, \ p_{3}=2b \ (a\leqslant b) \\
    & \Rightarrow (2a-4)(2b-4)<16\\
    &\Rightarrow(a-2)(b-2)<4
\end{align*}Thus,\begin{align*}
    &a-2=1, \ b-2 = 3 \Rightarrow a=3, \ b=5 \Rightarrow   p_{2}=6, \ p_{3}=10  \\
    &a-2=1, \ b-2 = 2 \Rightarrow a=3, \ b=4 \Rightarrow   p_{2}=6, \ p_{3}=8    \\
    &a-2=1, \ b-2 = 1 \Rightarrow a=3, \ b=3 \Rightarrow   p_{2}=6, \ p_{3}=6  \\
    &a-2=0, \ b-2 = n \Rightarrow a=2, \ b=2+n \Rightarrow   p_{2}=4, \ p_{3}=2(2+n)   \\  
\end{align*}and we conclude

\begin{align*}
   \fbox{$ \begin{array}{lll}
    &   (p_{1},p_{2},p_{3}) = (4.6.10)\cdots\rm \ \mathbf{great \ rhombicosidodecahedron} \\
     &   (p_{1},p_{2},p_{3}) = (4.6.8)\cdots\rm \ \mathbf{great \ rhombicuboctahedron} \\
      &   (p_{1},p_{2},p_{3}) = 4.6^{2}\cdots\rm \ \mathbf{truncated \  octahedron }\\
       &   (p_{1},p_{2},p_{3}) =4^{3} \cdots\rm \ \mathbf{cube} \\
        &   (p_{1},p_{2},p_{3}) = 4^{2}.m\  (m\geqslant 4)\ \cdots\rm \ \mathbf{prism} \end{array}$}
\end{align*}

We note that this subcase covers precisely the polyhedra with \textbf{\emph{bipartite}} graphs, i.e., if \textbf{V} is the set of vertices of the polyhedron and if $\mathbf{V}=\mathbf{V_{1}}\bigcup \mathbf{V_{2}}$ while $\mathbf{V_{1}}\bigcap \mathbf{V_{2}}=\emptyset$ and each edge of the graph goes from $\mathbf{V_{1}}$ to $\mathbf{V_{2}}.$  Equivalently, each $p_{k}$ is \textbf{\emph{even}}.\footnote{We thank \textsc{Michael Josephy} for this observation.}
\subsubsection{All faces have at least five sides and one exactly five sides: $p_{1}=5\leqslant p_{2}\leqslant p_{3}$}

This is quite similar the the previous section.  Since\begin{align*}
    &5=p_{1}\leqslant p_{2}\leqslant p_{3}   \\
    &\Rightarrow \frac{1}{p_{2}}+\frac{1}{p_{3}}-\frac{3}{10}>0\\
    &\Rightarrow (3p_{2}-10)(3p_{3}-10)<100. 
\end{align*}Again, a configuration argument shows that \begin{align*}
    &      p_{2}=p_{3}                   \\
    &\Rightarrow (3p_{2}-2)^{2}<100\\
    &\Rightarrow 15\leqslant 3p_{2}<20\\
    &\Rightarrow p_{2}=5, \ 6,  
\end{align*}which gives
\begin{align*}
   \fbox{$ \begin{array}{lll}
    &   (p_{1},p_{2},p_{3}) = 5^{3}\cdots\rm \ \mathbf{regular \ dodecahedron} \\
     &   (p_{1},p_{2},p_{3}) = 5.6^{2}\cdots\rm \ \mathbf{truncated \ icosahedron} \\
       \end{array}$}
\end{align*}
\section{Summary of our results}

We present a list of the polyhedra we have found and the section of the proof where they were determined.

For the regular polyhedra:
\large
\begin{center}
\fbox{\textbf{Regular Polyhedra}}
\end{center}
\normalsize
\begin{center}
\begin{tabular}{|c|c|}
  \hline
 \textbf {NAME } &\textbf{ Section where found}\\ \hline
 \textbf{Tetrahedron} & 4.3.1 \\ \hline
 \textbf{Octahedron} &4.2.1 \\ \hline
 \textbf{Icosahedron} & 4.1.2\\ \hline
 \textbf{Cube} & 4.3.2 \\ \hline
 \textbf{Dodecahedron} &4.3.3 \\ \hline
\end{tabular}
\end{center}

For the Archimedean polyhedra:
\large
\begin{center}
\fbox{\textbf{Archimedean Polyhedra}}
\end{center}
\normalsize
\begin{center}
\begin{tabular}{|c|c|}
  \hline
  \textbf{NAME} &\textbf{Section where found}  \\
 
 \hline
\textbf{cuboctahedron} &4.2.1\\
 \hline
 \textbf{great rhombicosidodecahedron} &4.3.2\\
 \hline
 \textbf{great rhombicuboctahedron} &4.3.2\\
 \hline
 \textbf{icosidodecahedron} &4.2.1\\
 \hline
 \textbf{small rhombicosidodecahedron} &4.2.1\\
 \hline
 \textbf{small rhombicuboctahedron} &4.2.1\\
 \hline
 \textbf{snub cube}&4.1.2\\
 \hline
 \textbf{snub dodecahedron} &4.1.2\\
 \hline
 \textbf{truncated cube} &4.3.1\\
 \hline
 \textbf{truncated dodecahedron} &4.3.1\\
 \hline
 \textbf{truncated icosahedron} &4.3.3\\
 \hline
 \textbf{truncated octahedron} &4.3.2\\
 \hline
 \textbf{truncated tetrahedron} &4.3.1\\
 \hline
 \textbf{prisms} &4.3.1, \ 4.3.2\\
 \hline
 \textbf{antiprisms} &4.2.1\\
 \hline
\end{tabular}
\end{center}

And we have completed the topological proof of \textsc{Archimedes}' theorem.

We have \emph{not} demonstrated that the polyhedra enumerated in \textsc{Archimedes}' theorem are in fact \emph{constructable}.  Again, this is done in the works of \textsc{Cromwell} \cite{Crom} and \textsc{Lines} \cite{Lines}.
\section{Final remarks}

As in the case of the topological proof that there are five regular polyhedra, we have proven much more!  We have found all \textbf{\emph{semiregular maps}} on any homeomorph of the sphere, a result of great generality. Although the metric proofs are of great interest, intrinsically and historically, the topological proof shows that they appeal to unessential properties of their metric realizations and that, at the root of it all, \textsc{Archimedes}' theorem is a consequence of certain combinatorial relations among the numbers of vertices, edges, and faces.

One wonders what \textsc{Archimedes} would have thought of our proof of his theorem.  We hope that he would have liked it.
\subsubsection*{Acknowledgment}
Support from the Vicerrector\'{\i}a de Investigaci\'on of the 
University of Costa Rica is acknowledged.

\end{document}